# Right, left and double division on semigroups that are semilattices of groups

**R. A. R. Monzo**

**Abstract.** The binary products of right, left or double division on semigroups that are semilattices of groups give interesting groupoid structures that are in one-one correspondence with semigroups that are semilattices of groups. This work is inspired by the known one-one correspondence between groups and Ward quasigroups.

## 1. Introduction.

It appears in the literature that in 1930 M. Ward was the first to find a set of axioms on $(S,*)$ (a set $S$ with a binary operation $*$ , called here a *groupoid* ) that ensure the existence of a group binary operation $\circ$ on $S$ such that $x*y = x\circ y^{-1}$ [14]. Such a groupoid was called a "division groupoid" by Polonijo [12] and it is clear that division groupoids are quasigroups.

Over the next 63 years many other sets of axioms on a groupoid were found that make it a division groupoid, now commonly known as a Ward quasigroup ([1], [2], [4], [6], [7], [8], [11], [12], [13], [15]). Perhaps the most impressive of these characterisations of Ward quasigroups is that of Higman and Neumann who found a single law making a groupoid a Ward quasigroup [7]. It is now known that a quasigroup is a Ward quasigroup if and only if it satisfies the law of right transitivity, $(x*z)*(y*z) = x*y$ [11]. It follows that a quasigroup is the dual of a Ward quasigroup, which we will call a Ward dual quasigroup, if and only if it satisfies the identity $(z*x)*(z*y) = x*y$.

Starting from any group $(G,\circ)$ we can form a Ward Quasigroup $(G,*)$ by defining $x*y = x\circ y^{-1}$; that is, $*$ is the operation of right division in the group $(G,\circ)$. Conversely, any Ward quasigroup $(W,*)$ is unipotent and its only idempotent $e = e*e = x*x$ (for any $x\in W$), is a right identity element. If we then define $(W,\circ)$ as $x\circ y = x*(e*y)$, $(W,\circ)$ is a group, $x^{-1} = e*x$ and $x*y = x\circ y^{-1}$. These mappings, $(G,\circ)\mapsto(G,*)$ and $(W,*)\mapsto(W,\circ)$ are inverse mappings, which implies that groups are in one-to-one correspondence with Ward quasigroups. (All this is well known.) In addition, a Ward quasigroup is an inverse groupoid, with the unique inverse of $x$ being $x^{-1} = e*x$. That is, the inverse of an element of a Ward quasigroup is the inverse element in the group it "induces".

In 2007 N. C. Fiala proved that a quasigroup $(S,*)$ satisfies the identity $[(e*e)*(x*z)]*[(e*y)*z]=x*y$ (for some $e\in S$ ) if and only if there is a group $(S,\circ)$ with identity element *e* such that $x*y=x^{-1}\circ y^{-1}$ [5]. Fiala called such groupoids "double Ward quasigroups". He noted that the binary operation $\circ$ on a double Ward quasigroup $S$ defined by $x\circ y=(e*x)*(e*y)$ is a group operation and that double Ward quasigroups are in one-to-one correspondence with groups. Double Ward quasigroups are also inverse quasigroups, with $x^{-1}=x$.

Our intention here is to explore the operations $x*y=x\bullet y^{-1}$ (called "right division"), $x*y=x^{-1}\bullet y$ (called "left division") and $x*y=x^{-1}\bullet y^{-1}$ (called "double division") when $(S,\bullet)$ is a semigroup and a semilattice of groups, where $x^{-1}$ is the inverse of $x$ in the group to which it belongs. We will prove that, for each of these operations, the collection of all such structures are in one-one correspondence with the collection of all semigroups that are semilattices of groups. A groupoid $(S,*)$ is in one of these collections if and only if there is a semigroup semilattice of groups for which the binary operation $*$ is right, left or double division. In these two senses, our results "extend" the results about Ward and double Ward quasigroups.

## 2. Preliminary definitions and known results

We call the pair $(S,*)$ a *groupoid* if $S$ is a set and $*$ is a binary operation on $S$ . We call $(T,*)$ a sub-groupoid of $(S,*)$ if $T\subseteq S$ and $T*T\subseteq T$ . We call an element $x$ of a groupoid $(S,*)$ *idempotent* if $x*x=x$. Then $E(S,*)=\{x\in S : x=x*x\}$ is the set of idempotent elements of $(S,*)$, which is not necessarily a sub-groupoid of $(S,*)$. The groupoid $(S,*)$ is called an *idempotent groupoid* [*a semilattice groupoid*] if all of its elements are idempotent [idempotent and commute]. A semilattice groupoid $(S,*)$ is called a *semigroup semilattice* if $(S,*)$ is a semigroup. A groupoid $(S,*)$ is called an (*idempotent*) *groupoid* $(T,\bullet)$ *of groupoids* $(T_\alpha,*/_{T_\alpha})$ $(\alpha\in T)$ if $S$ is a disjoint union of the $T_\alpha$ $(\alpha\in T)$ and $T_\alpha*T_\beta\subseteq T_{\alpha\bullet\beta}$ for all $\alpha,\beta\in T$. Note that $\bullet$ or $*$ need not be associative. We write $x_\alpha$ for an element of $T_\alpha$.

We call the groupoid $(S,*)$ *right* [*left*] *solvable* if for any $a,b\in S$ there exists $x\in S$ such that $a*x=b$ $[x*a=b]$. The groupoid $(S,*)$ is a *quasigroup* if it is right and left solvable, in which case it is right and left cancellative. We call the quasigroup $(S,*)$ a *Ward quasigroup* [*Ward dual quasigroup*] if it satisfies the identity $(x*z)*(y*z)=x*y$ $[(z*x)*(z*y)=x*y]$. The quasigroup $(S,*)$ is called a *double Ward quasigroup* if it satisfies the identity $\left[(e*e)*(x*z)\right]*\left[(e*y)*z\right]=x*y$ for some $e\in S$. A groupoid $(S,*)$ is called *an inverse groupoid* if for all $x\in S$ there exists a unique element $x^{-1}\in S$ such that $\left(x*x^{-1}\right)*x=x$ and $\left(x^{-1}*x\right)*x^{-1}=x^{-1}$. We let $(S,*)\cong(T,\bullet)$ denote that $(S,*)$ and $(T,\bullet)$ are isomorphic groupoids. If $(S,*)$ is a groupoid and $T\subseteq S$ then $*/_T:T\times T\to T$, where $s*/_T t=s*t$ $(s,t\in T)$. The *dual mapping* $\overleftarrow{\bullet}$ of a mapping $\bullet$ is defined as $x\overleftarrow{\bullet}y=y\bullet x$. The collection $\overleftarrow{\mathcal{C}}$ is the collection of all groupoids $\left(S,\overleftarrow{\bullet}\right)$, where $(S,\bullet)\in\mathcal{C}$. Clearly, $\mathcal{C}$ is in one-one correspondence with $\overleftarrow{\mathcal{C}}$.

A Ward quasigroup $(S,*)$ satisfies the following identities ([2], [9], [11]):

$$(x*z)*(y*z)=x*y \tag{1}$$

$$x*x=y*y=r \tag{2}$$

$$x*r=x \tag{3}$$

$$r*(x*y)=y*x \tag{4}$$

$$r*(r*x)=x \tag{5}$$

$$(x*y)*z=x*\left[z*(r*y)\right] \tag{6}$$

Dually, it follows that a Ward dual quasigroup $(S,\bullet)$ satisfies the following identities:

$$(z*x)*(z*y)=x*y \tag{7}$$

$$x*x=y*y=r \tag{8}$$

$$r*x=x \tag{9}$$

$$x*y=(y*x)*r \tag{10}$$

$$(x*r)*r=x \tag{11}$$

$$x*(y*z)=\left[(y*r)*x\right]*z \tag{12}$$

From Result 3 below, a double Ward quasigroup $(S,*)$ satisfies the following identities:

$$\left[(e*e)*(x*z)\right]*\left[(e*y)*z\right]=x*y \tag{13}$$

$$e*e=e \tag{14}$$

$$\left[e*(x*z)\right]*\left[(e*y)*z\right]=x*y \tag{15}$$

$$(y*x)*y=x=y*(x*y) \tag{16}$$

$$e*x=x*e \tag{17}$$

$$x*(x*e)=e=(e*x)*x \tag{18}$$

$$x*y=e*\left[(e*y)*(e*x)\right] \tag{19}$$

Note that a Ward [Ward dual] quasigroup $(S,*)$ has a unique right [left] identity element $r$. So, we will denote this by $(W,*,r)$ [$(WD,*,r)$]. We will denote a double Ward quasigroup by $(DW,*,e)$, although we note that the element $e$ may not be unique. The following results are well-known, or follow readily from well-known results (cf. eg. [2], [5], [11]).

**Result 1.** *If* $(W,*,r)$ *is a Ward quasigroup if and only if there is a group* $(W,\circ,r)$ *such that* $x*y=x\circ y^{-1}$ *for all* $x,y\in W$.

**Result 2.** *If* $(WD,*,r)$ *is a Ward dual quasigroup if and only if there is a group* $(WD,\circ,r)$ *such that* $x*y=y\circ x^{-1}$ *for all* $x,y\in W$.

**Result 3.** *If* $(DW,*,e)$ *is a double Ward quasigroup if and only if there is a group* $(DW,\circ,e)$ *such that* $x*y=x^{-1}\circ y^{-1}$ *for all* $x,y\in W$.

**Result 4.** *If* $(W,*,r)$ *is a Ward quasigroup then* $(W,\circ)$ *defined as* $x\circ y=x*(r*y)$ $(x,y\in W)$ *is a group with identity* $r$ *and* $x^{-1}=r*x$.

**Result 5.** *If* $(WD,*,r)$ *is a Ward dual quasigroup then* $(WD,\circ)$ *defined as* $x\circ y=(x*r)*y$ $(x,y\in W)$ *is a group with identity* $r$ *and* $x^{-1}=x*r$.

**Result 6.** *If* $(DW,*,e)$ *is a double Ward quasigroup then* $(DW,\circ)$ *defined as* $x\circ y=(e*x)*(e*y)$ $(x,y\in W)$ *is a group with identity* $e$ *and* $x^{-1}=e*x$.

Proof. (This result was noted in [5] without proof.) By definition we have $(x\circ y)\circ z=$ $\{e*[(e*x)*(e*y)]\}*(e*z)\overset{19}{=}(y*x)*(e*z)$ and $x\circ(y\circ z)=(e*x)*\{e*[(e*y)*(e*z)]\}\overset{19}{=}(e*x)*(z*y)$. Since, by 16 and 17, $x=(y*x)*y$ and $z=e*(e*z)$, we have $x\circ(y\circ z)=(e*x)*(z*y)=\{e*[(y*x)*y]\}*\{[e*(e*z)]*y\}\overset{15}{=}$ $\overset{15}{=}(y*x)*(e*z)=(x\circ y)\circ z$. So, $(DW,*)$ is a semigroup.

Now suppose that $a,b\in DW$. Since $(DW,*,e)$ is a quasigroup, there exists a unique $x\in DW$ such that $x*a=e*b$. So $a\circ x=(e*a)*(e*x)\overset{16,17,19}{=}e*(x*a)=e*(e*b)\overset{16,17}{=}b$. If $a\circ y=b$ then $b=(e*a)*(e*y)\overset{16,17,19}{=}e*(y*a)$ and so, by 16 and 17, $e*b=b*e=[e*(y*a)]*e=y*a$. But $x$ was unique, so $x=y$. Similarly, there exists a unique element $z\in DW$ such that $z\circ a=b$. So, $(DW,\circ)$ is an associative quasigroup; hence, it is a group. The facts that $e$ is the identity and $x^{-1}=e*x$ follow from 16 and 18. ■

**Definition.** In Results 1 – 6 we say that the quasigroup $(Q,*)$ in question is *induced by the group* $(Q,\circ)$ and that $(Q,\circ)$ *generates* $(Q,*)$.

As a consequence of Results 1 through 6, we have the following Corollaries:

**Corollary 7 ([2], [11]).** *The collection of all Ward quasigroups is in one-to-one correspondence with the collection of all groups.*

**Corollary 8.** *The collection of all Ward dual quasigroups is in one-to-one correspondence with the collection of all groups.*

**Corollary 9 [5].** *The collection of all double Ward quasigroups is in one-to-one correspondence with the collection of all groups.*

The following two results follow readily from Results 1 and 3 and proofs are omitted.

**Result 10.** *If* $(W,*,r)$ *is a Ward quasigroup and if we define* $x\bullet y=(r*x)*y$ *then* $(W,\bullet,r)$ *is a double Ward quasigroup.*

**Result 11.** *If* $(DW,*,e)$ *is a double Ward quasigroup and if we define* $x\bullet y=(e*x)*y$ *then* $(DW,\bullet,e)$ *is a Ward quasigroup.*

**Result 12.** *If* $(S,*)$ *is a semigroup semilattice* $Y$ *of Ward quasigroups* $\left(W_\alpha,*/_{W_\alpha},e_\alpha\right)$ $(\alpha\in Y)$ *and satisfies the identity* $(x_\alpha*y_\beta)*(z_\sigma*w_\gamma)=\left[x_\alpha*\left(w_\gamma^{-1}*y_\beta^{-1}\right)\right]*z_\sigma$ *then* $E(S,*)=\{e_\alpha:\alpha\in Y\}$, $e_\alpha*e_\beta=e_{\alpha\beta}$ *and the mapping* $\Psi:\left(\{e_\alpha:\alpha\in Y\},*/_{\{e_\alpha:\alpha\in Y\}}\right)\to Y$ *defined by* $\Psi e_\alpha=\alpha$ *is an isomorphism* (*that is,* $\left(E(S,*),*/_{E(S,*)}\right)\cong Y$ ).

Proof. It is clear from (2) that $E(S,*)=\{e_\alpha:\alpha\in Y\}$. We note that each $\left(W_\alpha,*/_{W_\alpha},e_\alpha\right)$ is an inverse groupoid, with $x_\alpha^{-1}=e_\alpha*x_\alpha$. Since a semigroup semilattice groupoid of inverse groupoids is an inverse groupoid, $(S,*)$ is an inverse groupoid. Hence, the identity $(x_\alpha*y_\beta)*(z_\sigma*w_\gamma)=\left[x_\alpha*\left(w_\gamma^{-1}*y_\beta^{-1}\right)\right]*z_\sigma$ has a clear meaning. We call this identity "(*I*)". Now, by definition, $e_\alpha*e_\beta\in W_{\alpha\beta}$. Therefore, $e_{\alpha\beta}\overset{(2)}{=}(e_\alpha*e_\beta)*(e_\alpha*e_\beta)\overset{(I)}{=}(e_\alpha*e_\beta)*e_\alpha$. Then,

$(e_\alpha*e_\beta)*e_\beta\overset{(I)}{=}\left[e_\alpha*(e_\beta*e_\alpha)\right]*e_\beta\overset{(I)}{=}(e_\alpha*e_\alpha)*(e_\beta*e_\beta)=(e_\alpha*e_\beta)=(e_\alpha*e_\beta)*e_{\alpha\beta}=(e_\alpha*e_\beta)*\left[(e_\alpha*e_\beta)*e_\alpha\right]\overset{(I)}{=}$

$\overset{(I)}{=}\left[e_\alpha*\left(e_\alpha*e_\beta\right)\right]*\left(e_\alpha*e_\beta\right)\overset{(I)}{=}\left\{\left[e_\alpha*\left(e_\beta*e_\alpha\right)\right]*e_\alpha\right\}*\left(e_\alpha*e_\beta\right)\overset{(I)}{=}\left[\left(e_\alpha*e_\beta\right)*e_\alpha\right]*\left(e_\alpha*e_\beta\right)=e_{\alpha\beta}*\left(e_\alpha*e_\beta\right).$

So,

$\left(e_\alpha*e_\beta\right)*e_\beta=\left(e_\alpha*e_\beta\right)=e_{\alpha\beta}*\left(e_\alpha*e_\beta\right)=\left[\left(e_\alpha*e_\beta\right)*e_\alpha\right]*\left[\left(e_\alpha*e_\beta\right)*e_\beta\right]\overset{(I)}{=}\left[\left(e_\alpha*e_\beta\right)*\left(e_\beta*e_\alpha\right)\right]*\left(e_\alpha*e_\beta\right).$

But, since $\left[\left(e_\alpha*e_\beta\right)*\left(e_\beta*e_\alpha\right)\right]\in W_{\alpha\beta}$, $W_{\alpha\beta}$ is a Ward quasigroup and $e_{\alpha\beta}*\left(e_\alpha*e_\beta\right)=$ $\left[\left(e_\alpha*e_\beta\right)*\left(e_\beta*e_\alpha\right)\right]*\left(e_\alpha*e_\beta\right)$, $e_{\alpha\beta}=\left(e_\alpha*e_\beta\right)*\left(e_\beta*e_\alpha\right)\overset{(2)}{=}\left(e_\alpha*e_\beta\right)*\left(e_\alpha*e_\beta\right)$ and $\left(e_\alpha*e_\beta\right)=\left(e_\beta*e_\alpha\right)$. But this implies $e_\alpha*\left(e_\beta*e_\sigma\right)=e_\alpha*\left(e_\sigma*e_\beta\right)=\left(e_\alpha*e_\alpha\right)*\left(e_\sigma*e_\beta\right)\overset{(I)}{=}\left[e_\alpha*\left(e_\beta*e_\alpha\right)\right]*e_\sigma\overset{(I)}{=}\left(e_\alpha*e_\beta\right)*e_\sigma$. Then, $e_{\alpha\beta}=\left(e_\alpha*e_\beta\right)*\left(e_\alpha*e_\beta\right)\overset{(I)}{=}\left(e_\alpha*e_\beta\right)*e_\alpha=e_\alpha*\left(e_\beta*e_\alpha\right)=e_\alpha*\left(e_\alpha*e_\beta\right)=\left(e_\alpha*e_\alpha\right)*e_\beta=e_\alpha*e_\beta$. It follows that the mapping $e_\alpha\mapsto\alpha$ is an isomorphism between $E(S,*)$ and $Y$. ■

Dually, we have the following result:

**Result 13.** *If* $(S,*)$ *is a semigroup semilattice groupoid* $Y$ *of Ward dual quasigroups* $\left(W_\alpha,*/_{W_\alpha},e_\alpha\right)$ $(\alpha\in Y)$ *and satisfies the identity* $(x*y)*(z*w)=y*\left[\left(z^{-1}*x^{-1}\right)*w\right]$ *then* $E(S,*)=\{e_\alpha:\alpha\in Y\}$, $e_\alpha*e_\beta=e_{\alpha\beta}$ *and the mapping* $\Psi:\left(\{e_\alpha:\alpha\in Y\},*/_{\{e_\alpha:\alpha\in Y\}}\right)\to Y$ *defined by* $\Psi e_\alpha=\alpha$ *is an isomorphism* (*that is,* $\left(E(S,*),*/_{E(S,*)}\right)\cong Y$ ).

**Result 14.** *If* $(S,*)$ *is a semigroup semilattice groupoid* $Y$ *of double Ward quasigroups* $\left(DW_\alpha,*/_{DW_\alpha},e_\alpha\right)$ $(\alpha\in Y)$ *then*

(1) $\left(\{e_\alpha:\alpha\in Y\},*/_{\{e_\alpha:\,\alpha\in Y\}}\right)\cong Y$ *if and only if*

(2) *for all* $\alpha,\beta,\sigma,\gamma\in Y$, $\left(e_\alpha*e_\beta\right)*\left(e_\sigma*e_\gamma\right)=e_\beta*\left[\left(e_\gamma*e_\sigma\right)*e_\alpha\right]$ *if and only if*

(3) *the mapping* $e_\alpha\mapsto\alpha$ *is an isomorphism from* $\left(\{e_\alpha:\alpha\in Y\},*/_{\{e_\alpha:\,\alpha\in Y\}}\right)$ *to* $Y$.

Proof. $(1\Rightarrow 2)$ Suppose that $\Psi:\left(\{e_\alpha:\alpha\in Y\},*/_{\{e_\alpha:\,\alpha\in Y\}}\right)\to Y$ is an isomorphism. Then $\Psi\left[\left(e_\alpha*e_\beta\right)*\left(e_\sigma*e_\gamma\right)\right]=(\Psi e_\alpha)(\Psi e_\beta)(\Psi e_\sigma)(\Psi e_\gamma)\overset{Y\text{ a semigroup semilattice}}{=}(\Psi e_\beta)\left\{\left[\left((\Psi e_\gamma)\right)\left((\Psi e_\sigma)\right)\right](\Psi e_\alpha)\right\}=$ $\Psi\left\{e_\beta*\left[\left(e_\gamma*e_\sigma\right)*e_\alpha\right]\right\}$. Since $\Psi$ is one-one, $\left(e_\alpha*e_\beta\right)*\left(e_\sigma*e_\gamma\right)=e_\beta*\left[\left(e_\gamma*e_\sigma\right)*e_\alpha\right]$.

$(2\Rightarrow 3)$ First, we prove that $\left(e_\alpha*e_\beta\right)=\left(e_\beta*e_\alpha\right)$ for all $\alpha,\beta\in Y$. By hypothesis, we have

$$\left(e_\alpha*e_\beta\right)*\left(e_\sigma*e_\gamma\right)=e_\beta*\left[\left(e_\gamma*e_\sigma\right)*e_\alpha\right] \qquad (20).$$

This implies $\left(e_\alpha*e_\beta\right)*e_\gamma=e_\beta*\left(e_\gamma*e_\alpha\right)$ (21) and so

$$\left(e_\beta * e_\alpha\right) * e_\beta \overset{(21)}{=} e_\alpha * e_\beta \overset{(21)}{=} e_\alpha * \left(e_\beta * e_\alpha\right) \qquad (22)$$

Now, $e_{\alpha\beta} \overset{(16)}{=} \left[\left(e_\alpha * e_\beta\right) * e_{\alpha\beta}\right] * \left(e_\alpha * e_\beta\right) \overset{(21)}{=} \left[e_\beta * \left(e_{\alpha\beta} * e_\alpha\right)\right] * \left(e_\alpha * e_\beta\right) \overset{(20)}{=} \left\{e_\beta * \left[\left(e_\alpha * e_{\alpha\beta}\right) * e_\beta\right]\right\} * \left(e_\alpha * e_\beta\right) \overset{(21)}{=}$

$= \left\{e_\beta * \left[e_{\alpha\beta} * \left(e_\beta * e_\alpha\right)\right]\right\} * \left(e_\alpha * e_\beta\right) \overset{(20)}{=} \left\{e_\beta * \left\{\left[e_{\alpha\beta} * \left(e_\alpha * e_\beta\right)\right] * e_{\alpha\beta}\right\}\right\} * \left(e_\alpha * e_\beta\right) \overset{(16)}{=} \left[e_\beta * \left(e_\alpha * e_\beta\right)\right] * \left(e_\alpha * e_\beta\right) \overset{(22)}{=}$

$= \left(e_\beta * e_\alpha\right) * \left(e_\alpha * e_\beta\right) \overset{(20)}{=} e_\alpha * \left[\left(e_\beta * e_\alpha\right) * e_\beta\right] \overset{(22)}{=} e_\alpha * \left(e_\alpha * e_\beta\right) \overset{(20)}{=} e_\alpha * \left[\left(e_\beta * e_\alpha\right) * e_\alpha\right] \overset{(22)}{=} e_\alpha * \left[e_\alpha * \left(e_\alpha * e_\beta\right)\right] = e_\alpha * e_{\alpha\beta}$

so,

$$e_{\alpha\beta} = e_\alpha * e_{\alpha\beta} = e_\alpha * \left(e_\alpha * e_\beta\right) = \left(e_\beta * e_\alpha\right) * \left(e_\alpha * e_\beta\right) \qquad (23)$$

Then, $$e_{\alpha\beta} * e_\alpha \overset{(22)}{=} e_{\alpha\beta} * \left(e_\alpha * e_{\alpha\beta}\right) \overset{(23)}{=} e_{\alpha\beta} * e_{\alpha\beta} \overset{(14)}{=} e_{\alpha\beta} \overset{(23)}{=} e_\alpha * \left(e_\alpha * e_\beta\right) \qquad (24)$$

Since, by (23), $e_{\alpha\beta} = \left(e_\beta * e_\alpha\right) * \left(e_\alpha * e_\beta\right)$, it follows from (16) that

$$\left(e_\alpha * e_\beta\right) * e_{\alpha\beta} = \left(e_\beta * e_\alpha\right) \qquad (25)$$

So, $\left(e_\beta * e_\alpha\right) = \left(e_\alpha * e_\beta\right) * e_{\alpha\beta} \overset{(21)}{=} e_\beta * \left(e_{\alpha\beta} * e_\alpha\right) \overset{(24)}{=} e_\beta * e_{\alpha\beta} \overset{(23)}{=} e_\beta * \left[e_\alpha * \left(e_\alpha * e_\beta\right)\right] \overset{(21)}{=}$

$= e_\beta * \left[\left(e_\beta * e_\alpha\right) * e_\alpha\right] \overset{(20)}{=} \left(e_\alpha * e_\beta\right) * \left(e_\alpha * e_\beta\right)$, which means that $\left(e_\alpha * e_\beta\right) * \left(e_\alpha * e_\beta\right) = \left(e_\beta * e_\alpha\right) \overset{(25)}{=}$

$\overset{(25)}{=} \left(e_\alpha * e_\beta\right) * e_{\alpha\beta}$. Since $\left(DW_{\alpha\beta}, */_{DW_{\alpha\beta}}, e_{\alpha\beta}\right)$ is a quasigroup and $Y$ a semilattice $e_{\beta\alpha} = e_{\alpha\beta} = e_\alpha * e_\beta = e_\beta * e_\alpha$. Also, $\left(e_\alpha * e_\beta\right) * e_\gamma = e_{\alpha\beta} * e_\gamma = e_{(\alpha\beta)\gamma} = e_{\alpha(\beta\gamma)} = e_\alpha * \left(e_\beta * e_\gamma\right)$. So, $\left(\left\{e_\alpha : \alpha \in Y\right\}, *\right)$ is a semigroup semilattice. Finally, the mapping $\Psi : \left(\left\{e_\alpha : \alpha \in Y\right\}, *\right) \to Y$ defined as $\Psi e_\alpha = \alpha$ satisfies $\Psi\left(e_\alpha * e_\beta\right) = \Psi e_{\alpha\beta} = \alpha\beta = \left(\Psi e_\alpha\right)\left(\Psi e_\beta\right)$ and so, since $\Psi$ is clearly one-one and onto $Y$, $\Psi$ is an isomorphism.

$(3 \Rightarrow 1)$ This is obvious. ■

## 3. Structures in one-to-one correspondence with semigroup semilattices of groups

We have seen that Ward quasigroups, Ward dual quasigroups and double Ward quasigroups are in one-to-one correspondence with groups. In this section, we extend these results to semigroups that are semilattices of groups. Note that in semigroup theory a "semilattice", a "union of groups" and a "semilattice of groups" are, by definition, semigroups. However, the definition of a semilattice (or idempotent) groupoid $(T, \bullet)$ of groupoids $\left(T_\alpha, */_{T_\alpha}\right)$ $(\alpha \in T)$ results in structures that are not necessarily associative, even when the $T_\alpha$ $(\alpha \in T)$ are all groups. Therefore, we use the terms "semigroup semilattice", "semigroup union of groups" and "semigroup semilattice of groups", terms that are redundant for semigroup theorists.

The idea is a straightforward one. We simply "extend" the binary product that gives the bijection between groups and Ward quasigroups, for example, to the semigroup semilattice of groups and to the resultant structure(s). So, we are working with structures that result from defining binary operations on a semigroup semilattice of groups $(S,\bullet)$ as follows: $x*y=x\bullet y^{-1}$ (called "right division") , $x*y=x^{-1}\bullet y$ (called "left division") and $x*y=x^{-1}\bullet y^{-1}$ (called "double division"). This is possible because a semigroup semilattice of groups is an inverse semigroup; that is, each element $x\in S$ has a unique inverse $x^{-1}$ that is the inverse of the element $x$ in the group to which it belongs [3, Theorem 4.11].

On the resultant structures $(S,*)$ we define binary operations as follows, respectively: $x_\alpha\otimes y_\beta=x_\alpha*\left(e_{\alpha\beta}*y_\beta\right)$, $x_\alpha\otimes y_\beta=\left(x_\alpha*e_{\alpha\beta}\right)*y_\beta$ and $x_\alpha\otimes y_\beta=\left(e_{\alpha\beta}*x_\alpha\right)*\left(e_{\alpha\beta}*y_\beta\right)$. These structures $(S,\otimes)$ turn out to be semigroup semilattices of groups. In each of these three cases, the mappings $(S,\bullet)\mapsto(S,*)$ and $(S,*)\mapsto(S,\otimes)$ are inverse mappings. Hence, we find three different collections of structures, each of which is in one-to-one correspondence with the collection $\mathbb{SLG}$ of all semigroup semilattices of groups.

**Lemma 3.1 [3, Theorem 4.11].** *A semigroup $(S,\bullet)$ is a semigroup semilattice $Y$ of groups $\left(G_\alpha,\bullet/_{G_\alpha},e_\alpha\right)$ $(\alpha\in Y)$ if and only if $(S,\bullet)$ is a semigroup union of groups and has commuting idempotents if and only if $(S,\bullet)$ is an inverse semigroup that is a semigroup union of groups if and only if $(S,\bullet)$ is a semigroup and a semigroup semilattice $Y\cong E(S,\bullet)$ of groups.*

Note that the following identity holds in an inverse semigroup $(S,\bullet)$ :

$$(x\bullet y)^{-1}=y^{-1}\bullet x^{-1} \qquad (26)$$

If $(S,\bullet)$ is a semigroup and a semilattice $Y$ of groups then it follows from Lemma 3.1 that

$$e_\alpha\bullet e_\beta=e_{\alpha\beta}=e_{\beta\alpha}=e_\beta\bullet e_\alpha \quad \text{for all } (\alpha\in Y) \qquad (27)$$

**Lemma 3.2.** *Suppose that $(S,\bullet)$ is a semigroup semilattice $Y$ of groups $(G_\alpha,e_\alpha)(\alpha\in Y)$ and that $x_\alpha*y_\beta=x_\alpha\bullet y_\beta^{-1}$ for all $x_\alpha\in W_\alpha$, $y_\beta\in W_\beta$ and all $\alpha,\beta\in Y$. Then*

(3.21) $(S,*)$ *is an inverse groupoid with* $x_\alpha^{-1}=e_\alpha*x_\alpha$ $(\alpha\in Y)$ ,

(3.22) $E(S,*)\cong E(S,\bullet)\cong Y$,

(3.23) $(S,*)$ *is a semigroup semilattice $Y$ of Ward quasigroups* $\left(G_\alpha,*/_{G_\alpha},e_\alpha\right)(\alpha\in Y)$,

(3.24) $\left(x_\alpha*y_\beta\right)*\left(z_\sigma*w_\gamma\right)=\left[x_\alpha*\left(w_\gamma^{-1}*y_\beta^{-1}\right)\right]*z_\sigma$,

(3.25) $x_\alpha * (e_{\alpha\beta} * y_\beta) = x_\alpha * y_\beta^{-1}$ and

(3.26) $(x_\alpha * y_\beta) = (y_\beta * x_\alpha)^{-1}$

Proof. (3.21) It is straightforward to calculate that $x_\alpha^{-1}$, the inverse of $x_\alpha$ in the group to which it belongs, is the unique inverse of $x_\alpha$ in $(S,*)$. That is, $x_\alpha^{-1} = e_\alpha * x_\alpha$.

(3.22) $x_\alpha = x_\alpha * x_\alpha$ if and only if $x_\alpha = x_\alpha \bullet x_\alpha^{-1} = e_\alpha$, the identity of the group to which $x_\alpha$ belongs. Then, $e_\alpha * e_\beta = e_\alpha \bullet e_\beta^{-1} = e_\alpha \bullet e_\beta \overset{(21)}{=} e_{\alpha\beta}$. Since, by Lemma 3.1, in $(S,\bullet)$ we have $E(S,\bullet) \cong Y$, $E(S,*) \cong E(S,\bullet) \cong Y$.

(3.23) Since $x_\alpha * y_\beta = x_\alpha \bullet y_\beta^{-1} \in G_\alpha \bullet G_\beta \subseteq G_{\alpha\beta}$, $G_\alpha * G_\beta \subseteq G_{\alpha\beta}$. Since $x_\alpha * y_\alpha = x_\alpha \bullet y_\alpha^{-1}$ in each $(G_\alpha, */_{G_\alpha}, e_\alpha)$, by Result 1, $(G_\alpha, */_{G_\alpha}, e_\alpha)$ is a Ward quasigroup for all $\alpha \in Y$. By definition then, $(S,*)$ is a semigroup semilattice $Y$ of Ward quasigroups $(G_\alpha, */_{G_\alpha}, e_\alpha)$ $(\alpha \in Y)$.

(3.24) Using the facts that $x_\alpha * y_\beta = x_\alpha \bullet y_\beta^{-1}$ and $(x_\alpha \bullet y_\beta)^{-1} \overset{(26)}{=} y_\beta^{-1} \bullet x_\alpha^{-1}$ it is straightforward to calculate that $(x_\alpha * y_\beta) * (z_\sigma * w_\gamma) = x_\alpha \bullet y_\beta^{-1} \bullet w_\gamma \bullet z_\sigma^{-1} = [x_\alpha * (w_\gamma^{-1} * y_\beta^{-1})] * z_\sigma$.

(3.25) $x_\alpha * (e_{\alpha\beta} * y_\beta) = x_\alpha \bullet (e_{\alpha\beta} \bullet y_\beta^{-1})^{-1} \overset{(26)}{=} x_\alpha \bullet (y_\beta \bullet e_{\alpha\beta}^{-1}) = (x_\alpha \bullet y_\beta) \bullet e_{\alpha\beta} = x_\alpha \bullet y_\beta = x_\alpha * y_\beta^{-1}$.

(3.26) $(x_\alpha * y_\beta)^{-1} = (x_\alpha \bullet y_\beta^{-1})^{-1} \overset{(26)}{=} y_\beta \bullet x_\alpha^{-1} = y_\beta * x_\alpha$. ■

**Definition.** If $(S,\bullet)$ is a semigroup semilattice $Y$ of groups $(G_\alpha, e_\alpha)$ $(\alpha \in Y)$ and $x * y = x \bullet y^{-1}$ then we denote $(S,*)$ by $SLWQ(S,\bullet)$. We define $\mathbb{SLWQ}$ as the collection of all semigroup semilattices $Y$ of Ward quasigroups $(W_\alpha, */_{G_\alpha}, e_\alpha)$ $(\alpha \in Y)$ that satisfy 3.24. In particular, $SLWQ(S,\bullet) \in \mathbb{SLWQ}$. Note once again that a semigroup semilattice of inverse groupoids is an inverse groupoid. So, 3.24 (3.25 and 3.26) have a clear meaning.

**Lemma 3.3.** *Suppose that $(S,*)$ is a semigroup semilattice $Y$ of Ward quasigroups $(W_\alpha, */_{W_\alpha}, e_\alpha)$ $(\alpha \in Y)$ and satisfies* 3.24. *Define $x_\alpha \bullet y_\beta = x_\alpha * (e_{\alpha\beta} * y_\beta)$. Then $(S,\bullet)$ is a semigroup and a semigroup semilattice $Y$ of groups $(W_\alpha, \bullet/_{W_\alpha}, e_\alpha)$ $(\alpha \in Y)$ with $Y \cong E(S,\bullet) \cong E(S,*)$.*

Proof. As previously noted in the proof of Result 12, since each $\left(W_\alpha, */_{W_\alpha}, e_\alpha\right)$ is an inverse groupoid, with $x_\alpha^{-1} = e_\alpha * x_\alpha$ and since a semigroup semilattice of inverse groupoids is an inverse groupoid, $(S,*)$ is an inverse groupoid.

We prove that 3.24 implies 3.25. We have

$$x_\alpha * \left(e_{\alpha\beta} * y_\beta\right) = \left(x_\alpha * e_\alpha\right) * \left(e_{\alpha\beta} * y_\beta\right) \overset{3.24}{=} \left[x_\alpha * \left(y_\beta^{-1} * e_\alpha^{-1}\right)\right] * e_{\alpha\beta} = \left[x_\alpha * \left(y_\beta^{-1} * e_\alpha\right)\right] = \left(x_\alpha * e_\alpha\right) * \left(y_\beta^{-1} * e_\alpha\right) \overset{3.24}{=}$$
$$= \left[x_\alpha * \left(e_\alpha * e_\alpha\right)\right] * y_\beta^{-1} = \left(x_\alpha * e_\alpha\right) * y_\beta^{-1} = x_\alpha * y_\beta^{-1}$$ so, 3.25 is valid.

Next, we prove that 3.24 implies 3.26. Since $x_\alpha * y_\beta = x_\alpha * \left(y_\beta^{-1}\right)^{-1} \overset{3.25}{=} x_\alpha * \left(e_{\alpha\beta} * y_\beta^{-1}\right)$,

$$x_\alpha * y_\beta = x_\alpha * \left(e_{\alpha\beta} * y_\beta^{-1}\right) \overset{5}{=} \left[e_\alpha * \left(e_\alpha * x_\alpha\right)\right] * \left(e_{\alpha\beta} * y_\beta^{-1}\right) \overset{3.24}{=} \left[e_\alpha * \left(y_\beta * x_\alpha\right)\right] * e_{\alpha\beta} = \left[e_\alpha * \left(y_\beta * x_\alpha\right)\right] * \left(e_{\alpha\beta} * e_{\alpha\beta}\right) \overset{3.24}{=}$$
$$= \left\{e_\alpha * \left[e_{\alpha\beta} * \left(y_\beta * x_\alpha\right)^{-1}\right]\right\} * e_{\alpha\beta} = e_\alpha * \left(y_\beta * x_\alpha\right) = \left(e_\alpha * e_\alpha\right) * \left[\left(y_\beta * x_\alpha\right) * e_{\alpha\beta}\right] \overset{3.24}{=} \left[e_\alpha * \left(e_{\alpha\beta} * e_\alpha\right)\right] * \left(y_\beta * x_\alpha\right) \overset{\text{Result 12}}{=}$$
$$= e_{\alpha\beta} * \left(y_\beta * x_\alpha\right) = \left(y_\beta * x_\alpha\right)^{-1}$$ so, 3.26 is valid.

Now $x_\alpha = x_\alpha \bullet x_\alpha$ if and only if $x_\alpha = x_\alpha * \left(e_\alpha * x_\alpha\right) = x_\alpha * e_\alpha$ if and only if $e_\alpha = e_\alpha * x_\alpha = x_\alpha * x_\alpha$ if and only if $x_\alpha = e_\alpha$. Also, we have $e_\alpha * e_\beta = \left[e_\alpha * \left(e_\alpha * e_\alpha\right)\right] * e_\beta \overset{3.24}{=} e_\alpha * \left(e_\beta * e_\alpha\right)$. Then,

$$e_\alpha \bullet e_\beta = e_\alpha * \left(e_{\alpha\beta} * e_\beta\right) = \left(e_\alpha * e_\alpha\right) * \left(e_{\alpha\beta} * e_\beta\right) \overset{3.24}{=} \left[e_\alpha * \left(e_\beta * e_\alpha\right)\right] * e_{\alpha\beta} = e_\alpha * \left(e_\beta * e_\alpha\right) \overset{3.24}{=} e_\alpha * e_\beta$$ and so $\bullet/_{E(S,*)} = */_{E(S,*)}$ and $E(S,\bullet) \cong E(S,*)$. Using Result 12, $E(S,\bullet) \cong E(S,*) \cong Y$ is a semigroup semilattice. Since, for each $\left(W_\alpha, */_{W_\alpha}, e_\alpha\right)$, $x_\alpha \bullet y_\alpha = x_\alpha * \left(e_\alpha * y_\alpha\right)$, by Result 4, each $\left(W_\alpha, \bullet/_{W_\alpha}, e_\alpha\right)$ is a group. Since $x_\alpha \bullet y_\beta = x_\alpha * \left(e_{\alpha\beta} * y_\beta\right) \in W_{\alpha\beta}$, $W_\alpha \bullet W_\beta \subseteq W_{\alpha\beta}$, and so $(S,\bullet)$ is a semigroup semilattice $Y$ of groups. So, we only need to prove that $(S,\bullet)$ is a semigroup.

We have $\left(x_\alpha \bullet y_\beta\right) \bullet z_\gamma = \left[x_\alpha * \left(e_{\alpha\beta} * y_\beta\right)\right] * \left(e_{\alpha\beta\gamma} * z_\gamma\right) \overset{3.25}{=} \left(x_\alpha * y_\beta^{-1}\right) * \left(e_{\alpha\beta\gamma} * z_\gamma\right) \overset{3.24}{=}$
$= \left[x_\alpha * \left(z_\gamma^{-1} * y_\beta\right)\right] * e_{\alpha\beta\gamma} = \left[x_\alpha * \left(z_\gamma^{-1} * y_\beta\right)\right]$.

Then, $x_\alpha \bullet \left(y_\beta \bullet z_\gamma\right) = x_\alpha * \left\{e_{\alpha\beta\gamma} * \left[y_\beta * \left(e_{\beta\gamma} * z_\gamma\right)\right]\right\} \overset{3.25}{=} x_\alpha * \left[e_{\alpha\beta\gamma} * \left(y_\beta * z_\gamma^{-1}\right)\right] \overset{3.25}{=} x_\alpha * \left(y_\beta * z_\gamma^{-1}\right)^{-1}$
$\overset{3.26}{=} x_\alpha * \left(z_\gamma^{-1} * y_\beta\right) = \left(x_\alpha \bullet y_\beta\right) \bullet z_\gamma$ and so $(S,\bullet)$ is a semigroup. ■

**Corollary 3.4.** *If* $(S,*)$ *is a semigroup semilattice* $Y$ *of Ward quasigroups* $\left(W_\alpha, */_{G_\alpha}, e_\alpha\right)$ $(\alpha \in Y)$ *and satisfies* 3.24 *then*

(1) $(S,*)$ *satisfies* 3.25 *and* 3.26 *and*

(2) *there exists* $(S,\bullet) \in \mathbb{SLG}$ *such that* $x * y = x \bullet y^{-1}$ *for all* $x, y \in S$

Proof. Part (1) was proved in Lemma 3.3. For part (2), let $(S,\bullet)$ be as in Lemma 3.3. Then, as proved in 3.3, $(S,\bullet)\in \mathbb{SLG}$. Also, $x_\alpha \bullet y_\beta^{-1} = (x_\alpha * e_\alpha)*(e_{\alpha\beta} * y_\beta^{-1}) \overset{3.24}{=} \left[x_\alpha *(y_\beta * e_\alpha)\right]* e_{\alpha\beta} =$

$= \left[x_\alpha *(y_\beta * e_\alpha)\right] = \left[(x_\alpha * e_\alpha)*(y_\beta * e_\alpha)\right] \overset{3.24}{=} (x_\alpha * e_\alpha)* y_\beta = x_\alpha * y_\beta$. ■

**Definition.** If $(S,*)$ is a semigroup semilattice $Y$ of Ward quasigroups $\left(W_\alpha, */_{W_\alpha}, e_\alpha\right)(\alpha \in Y)$ and satisfies 3.24, and if we define $x_\alpha \bullet y_\beta = x_\alpha *(e_{\alpha\beta} * y_\beta)$ then we denote the semigroup semilattice $Y$ of groups $(S,\bullet)$ as $SLG(S,*)$.

**Theorem 3.5.** *For all* $(S,*)\in \mathbb{SLWQ}$, $SLWQ(SLG(S,*)) = (S,*)$ *and for all* $(S,\bullet)\in \mathbb{SLG}$, $SLG(SLWQ(S,\bullet)) = (S,\bullet)$.

Proof. In $SLG(S,*)$ the product is $x_\alpha \bullet y_\beta = x_\alpha *(e_{\alpha\beta} * y_\beta)$.The product $\otimes$ in $SWQ(SLG(S,*))$ is $x_\alpha \otimes y_\beta = x_\alpha \bullet y_\beta^{-1}$. So,

$x_\alpha * y_\beta = (x_\alpha * y_\beta)*(e_{\alpha\beta} * e_{\alpha\beta}) \overset{3.25}{=} \left[x_\alpha *(e_{\alpha\beta}^{-1} * y_\beta^{-1})\right]* e_{\alpha\beta} = x_\alpha *(e_{\alpha\beta} * y_\beta^{-1}) = x_\alpha \bullet y_\beta^{-1} = x_\alpha \otimes y_\beta$.

Therefore, $SWQ(SLG(S,*)) = (S,*)$.

In $SWQ(S,\bullet)$ the product is $x * y = x \bullet y^{-1}$. The product $\oplus$ in $SLG(SWQ(S,\bullet))$ is

$x_\alpha \oplus y_\beta = x_\alpha *(e_{\alpha\beta} * y_\beta) = x_\alpha \bullet (e_{\alpha\beta} \bullet y_\beta^{-1})^{-1} \overset{20}{=} x_\alpha \bullet y_\beta \bullet e_{\alpha\beta}^{-1} = x_\alpha \bullet y_\beta$. Hence, $SLG(SWQ(S,\bullet)) = (S,\bullet)$. ■

The second part of the following Corollary can be viewed as an "extension" of Result 1.

**Corollary 3.6.** *There is a one-to-one correspondence between semigroup semilattices of groups* ($\mathbb{SLG}$) *and groupoids* $(S,*)$ *that are semigroup semilattices* $Y$ *of Ward quasigroups* $\left(W_\alpha, */_{G_\alpha}, e_\alpha\right)(\alpha \in Y)$ *and that satisfy* 3.24 ($\mathbb{SLWQ}$). *Also,* $(S,*)\in \mathbb{SLWQ}$ *if and only if there exists* $(S,\bullet)\in \mathbb{SLG}$ *such that* $x * y = x \bullet y^{-1}$ *for all* $x, y \in S$

**Corollary 3.7.** *There is a one-to-one correspondence between semigroup semilattices of abelian groups and groupoids* $(S,*)$ *that are semigroup semilattices* $Y$ *of medial Ward quasigroups* $\left(W_\alpha, */_{G_\alpha}, e_\alpha\right)(\alpha \in Y)$ *and that satisfy* 3.24.

Proof. The proof here follows that of 3.2, 3.3 and 3.5, using the additional fact that a groupoid is a medial Ward quasigroup if and only if it is induced by an abelian group. ■

**Lemma 3.8.** *Suppose that* $(S,\bullet)$ *is a semigroup semilattice* $Y$ *of groups* $(G_\alpha, e_\alpha)\,(\alpha \in Y)$ *and that* $x * y = x^{-1} \bullet y$. *Then*

(3.81) $(S,*)$ *is an inverse groupoid with* $x_\alpha^{-1} = x_\alpha * e_\alpha \;\; (\alpha \in Y)$ ,

(3.82) $E(S,*) \cong E(S,\bullet) \cong Y$,

(3.83) $(S,*)$ *is a semigroup semilattice* $Y$ *of Ward dual quasigroups* $\left(G_\alpha, */_{G_\alpha}, e_\alpha\right)(\alpha \in Y)$,

(3.84) $(x_\alpha * y_\beta)*(z_\sigma * w_\gamma) = y_\beta * \left[\left(z_\sigma^{-1} * x_\alpha^{-1}\right) * w_\gamma\right]$,

(3.85) $(x_\alpha * e_{\alpha\beta}) * y_\beta = x_\alpha^{-1} * y_\beta$ and

(3.86) $(x_\alpha * y_\beta) = (y_\beta * x_\alpha)^{-1}$

Proof. Note that it follows from Lemma 3.1 that $\mathbb{SLG} = \overleftarrow{\mathbb{SLG}}$. Since $x \overleftarrow{*} y = y * x = y^{-1} \bullet x = x \overleftarrow{\bullet} y^{-1}$ and $\left(S, \overleftarrow{\bullet}\right) \in \mathbb{SLG}$, $\left(S, \overleftarrow{*}\right)$ satisfies 3.21 to 3.26. Hence, $(S,*)$ satisfies 3.81 to 3.86. ■

**Definition.** $SLWD(S,\bullet)$ will denote $(S,*)$ in Lemma 3.8 above. We denote the collection of all groupoids $(S,*)$ that are semigroup semilattices $Y$ of Ward dual quasigroups $\left(WD_\alpha, */_{WD_\alpha}, e_\alpha\right)(\alpha \in Y)$ that satisfy 3.84 as $\mathbb{SLWDQ}$. So, $SLWD(S,\bullet) = (S,*) \in \mathbb{SLWDQ}$.

**Lemma 3.9.** *Suppose that* $(S,*)$ *is a semigroup semilattice* $Y$ *of Ward dual quasigroups* $\left(WD_\alpha, */_{WD_\alpha}, e_\alpha\right)(\alpha \in Y)$ *and satisfies* 3.84. *Define* $x_\alpha \bullet y_\beta = (x_\alpha * e_{\alpha\beta}) * y_\beta$. *Then* $(S,\bullet)$ *is a semigroup and a semigroup semilattice* $Y$ *of groups.*

Proof. $\left(S, \overleftarrow{*}\right)$ is clearly a semigroup semilattice $Y$ of Ward quasigroups $\left(WD_\alpha, */_{WD_\alpha}, e_\alpha\right)(\alpha \in Y)$ and satisfies 3.24. Also, $x \overleftarrow{\bullet} y = (y * e_{\alpha\beta}) * x = x \overleftarrow{*} \left(e_{\alpha\beta} \overleftarrow{*} y\right)$. By Lemma 3.3, $\left(S, \overleftarrow{\bullet}\right)$ is a semigroup and a semigroup semilattice $Y$ of groups with $Y \cong E\left(S, \overleftarrow{\bullet}\right) \cong E\left(S, \overleftarrow{*}\right)$. Hence, $(S,\bullet)$ is a semigroup and a semigroup semilattice $Y$ of groups with $Y \cong E(S,\bullet) \cong E(S,*)$. ■

**Definition.** $SLG(S,*)$ will denote $(S,\bullet)$ in Lemma 3.9 above.

**Theorem 3.10.** *For all* $(S,\bullet)\in\mathbb{SLG}$, $SLG(SLWD(S,\bullet))=(S,\bullet)$ *and for all* $(S,*)\in\mathbb{SLWDQ}$, $SLWD(SLG(S,*))=(S,*)$.

Proof. The product in $SLWD(S,\bullet)$ is $x*y=x^{-1}\bullet y$. The product in $SLG(SLWD(S,\bullet))$ is $x_\alpha\otimes y_\beta=(x_\alpha*e_{\alpha\beta})*y_\beta$. So, $x_\alpha\otimes y_\beta=(x_\alpha*e_{\alpha\beta})*y_\beta=(x_\alpha{}^{-1}\bullet e_{\alpha\beta})*y_\beta\overset{20}{=}e_{\alpha\beta}\bullet x_\alpha\bullet y_\beta=x_\alpha\bullet y_\beta$. Hence, $SLG(SLWD(S,\bullet))=(S,\bullet)$.

The product in $SLG(S,*)$ is $x_\alpha\bullet y_\beta=(x_\alpha*e_{\alpha\beta})*y_\beta$. The product in $SLWD(SLG(S,*))$ is $x\oplus y=x^{-1}\bullet y=[(x_\alpha*e_\alpha)*e_{\alpha\beta}]*y_\beta\overset{3.84}{=}x_\alpha*y_\beta$ and so $SLWD(SLG(S,*))=(S,*)$. ■

**Corollary 3.11.** *There is a one-to-one correspondence between semigroup semilattices of groups* ($\mathbb{SLG}$) *and groupoids* $(S,*)$ *that are semigroup semilattices* $Y$ *of Ward dual quasigroups* $(WD_\alpha,*/_{WD_\alpha},e_\alpha)$ $(\alpha\in Y)$ *and that satisfy* 3.81, 3.83 and 3.84 ($\mathbb{SLWDQ}$). *Also,* $(S,*)\in\mathbb{SLWDQ}$ *if and only if there exists* $(S,\bullet)\in\mathbb{SLG}$ *such that* $x*y=x^{-1}\bullet y$ *for all* $x,y\in S$.

**Corollary 3.12.** *There is a one-to-one correspondence between semigroup semilattices of abelian groups and groupoids* $(S,*)$ *that are semigroup semilattices* $Y$ *of unipotent, left-unital right modular quasigroups* $(ULR_\alpha,*/_{WD_\alpha},e_\alpha)$ $(\alpha\in Y)$ *satisfying* 3.84.

Proof. A Ward dual quasigroup is a unipotent, left-unital right modular quasigroup if and only if it is medial if and only if it is induced by an abelian group. Using this fact, the proof of Corollary 3.12 exactly follows those of Lemmas 3.8, 3.9 and Theorem 3.10.

**Lemma 3.13.** *Suppose that* $(S,\bullet)$ *is a semigroup semilattice* $Y$ *of groups* $(G_\alpha,e_\alpha)$ $(\alpha\in Y)$ *and that* $x*y=x^{-1}\bullet y^{-1}$. *Then*

(3.131) $(\{e_\alpha:\alpha\in Y\},*)\cong(\{e_\alpha:\alpha\in Y\},\bullet)\cong Y$ *is a semigroup semilattice*

(3.132) $(S,*)$ *is a semigroup semilattice* $Y$ *of double Ward quasigroups* $(G_\alpha,*/_{G_\alpha},e_\alpha)$ $(\alpha\in Y)$,

(3.133) $(S,*)$ *satisfies the identity*

$\{e_{\alpha\beta\gamma}*[(e_{\alpha\beta}*x_\alpha)*(e_{\alpha\beta}*y_\beta)]\}*(e_{\alpha\beta\gamma}*z_\gamma)=(e_{\alpha\beta\gamma}*x_\alpha)*\{e_{\alpha\beta\gamma}*[(e_{\beta\gamma}*y_\beta)*(e_{\beta\gamma}*z_\gamma)]\}$ *and*

(3.134) $(S,*)$ *satisfies the identity* $[e_{\alpha\beta}*(e_\alpha*x_\alpha)]*[e_{\alpha\beta}*(e_\beta*y_\beta)]=x_\alpha*y_\beta$.

Proof. For any $e_\alpha, e_\beta \in Y$, $e_\alpha * e_\beta = e^{-1}{}_\alpha \bullet e_\beta{}^{-1} = e_\alpha \bullet e_\beta$. Then, $e_\alpha * e_\beta = e_\alpha \bullet e_\beta \overset{3.1}{=} e_\beta \bullet e_\alpha = e_{\alpha\beta}$. Hence, $\left(e_\alpha * e_\beta\right) * e_\sigma = e_{\alpha\beta} * e_\sigma = e_{(\alpha\beta)\sigma} = e_{\alpha(\beta\sigma)} = e_\alpha * \left(e_\beta * e_\gamma\right)$. By Lemma 3.1, $\left(\{e_\alpha : \alpha \in Y\}, *\right) \cong \left(\{e_\alpha : \alpha \in Y\}, \bullet\right) \cong Y$ is a semigroup semilattice and so 3.131 is valid.

Each $\left(G_\alpha, */_{G_\alpha}, e_\alpha\right)$ has product $x_\alpha * y_\alpha = x_\alpha{}^{-1} \bullet y_\alpha{}^{-1}$ and therefore, by Result 3, each $\left(G_\alpha, */_{G_\alpha}, e_\alpha\right)$ is a double Ward quasigroup. Since $x_\alpha * y_\beta = x_\alpha{}^{-1} \bullet y_\beta{}^{-1} \in G_{\alpha\beta}$, 3.132 is valid. Then, we have

$$\left\{e_{\alpha\beta\gamma} * \left[\left(e_{\alpha\beta} * x_\alpha\right) * \left(e_{\alpha\beta} * y_\beta\right)\right]\right\} * \left(e_{\alpha\beta\gamma} * z_\gamma\right) =$$

$$\left\{e_{\alpha\beta\gamma} * \left[\left(e_{\alpha\beta}{}^{-1} \bullet x_\alpha{}^{-1}\right) * \left(e_{\alpha\beta}{}^{-1} \bullet y_\beta{}^{-1}\right)\right]\right\} * \left(e_{\alpha\beta\gamma}{}^{-1} \bullet z_\gamma{}^{-1}\right) = \left\{e_{\alpha\beta\gamma} * \left(x_\alpha \bullet y_\beta\right)\right\} * \left(e_{\alpha\beta\gamma}{}^{-1} \bullet z_\gamma{}^{-1}\right) = \left[e_{\alpha\beta\lambda}{}^{-1} \bullet y_\beta{}^{-1} \bullet x_\alpha{}^{-1}\right]^{-1} \bullet z_\gamma \bullet e_{\alpha\beta\gamma} =$$

$$\left(x_\alpha \bullet y_\beta \bullet e_{\alpha\beta\gamma}\right) \bullet z_\gamma \bullet e_{\alpha\beta\gamma} = \left(x_\alpha \bullet y_\beta\right) \bullet \left[e_{\alpha\beta\gamma} \bullet \left(z_\gamma \bullet e_{\alpha\beta\gamma}\right)\right] = \left(x_\alpha \bullet y_\beta\right) \bullet \left(z_\gamma \bullet e_{\alpha\beta\gamma}\right) = \left(x_\alpha \bullet y_\beta \bullet z_\gamma\right) \bullet e_{\alpha\beta\gamma} = \left(x_\alpha \bullet y_\beta \bullet z_\gamma\right)$$

. Also,

$$\left(e_{\alpha\beta\gamma} * x_\alpha\right) * \left\{e_{\alpha\beta\gamma} * \left[\left(e_{\beta\gamma} * y_\beta\right) * \left(e_{\beta\gamma} * z_\gamma\right)\right]\right\} = \left(e_{\alpha\beta\lambda}{}^{-1} \bullet x_\alpha{}^{-1}\right)^{-1} * \left\{e_{\alpha\beta\gamma}{}^{-1} \bullet \left[\left(e_{\beta\gamma}{}^{-1} \bullet y_\beta{}^{-1}\right)^{-1} \bullet \left(e_{\beta\gamma}{}^{-1} \bullet z_\gamma{}^{-1}\right)^{-1}\right]^{-1}\right\}^{-1} =$$

$$= \left(x_\alpha \bullet e_{\alpha\beta\gamma}\right) \bullet e_{\alpha\beta\gamma} \bullet \left(y_\beta \bullet e_{\beta\gamma}\right) \bullet \left(z_\gamma \bullet e_{\beta\gamma}\right) = x_\alpha \bullet \left(e_\alpha \bullet e_{\beta\gamma}\right) \bullet \left\{y_\beta \bullet \left[e_{\beta\gamma} \bullet \left(z_\gamma \bullet e_{\beta\gamma}\right)\right]\right\} = \left(x_\alpha \bullet e_\alpha\right) \bullet \left\{e_{\beta\gamma}\left[\left(y_\beta \bullet z_\gamma\right) \bullet e_{\beta\gamma}\right]\right\} = x_\alpha \bullet y_\beta \bullet z_\gamma$$

and so 3.133 is valid. Finally,

$$\left[e_{\alpha\beta} * \left(e_\alpha * x_\alpha\right)\right] * \left[e_{\alpha\beta} * \left(e_\beta * y_\beta\right)\right] = \left[e_{\alpha\beta}{}^{-1} \bullet \left(e_\alpha{}^{-1} \bullet x_\alpha{}^{-1}\right)^{-1}\right] * \left[e_{\alpha\beta}{}^{-1} \bullet \left(e_\beta{}^{-1} \bullet y_\beta{}^{-1}\right)^{-1}\right] = e_\alpha \bullet x_\alpha{}^{-1} \bullet \left[e_{\alpha\beta} \bullet \left(e_\beta \bullet y_\beta{}^{-1} \bullet e_{\alpha\beta}\right)\right] =$$

$$= e_\alpha \bullet x_\alpha{}^{-1} \bullet e_\beta \bullet y_\beta{}^{-1} \bullet e_{\alpha\beta} = e_\alpha \bullet x_\alpha{}^{-1} \bullet e_\beta \bullet y_\beta{}^{-1} = x_\alpha{}^{-1} \bullet y_\beta{}^{-1} = x_\alpha * y_\beta$$

and 3.134 is valid. ■

**Definition.** $SLDWQ(S, \bullet)$ denotes $(S, *)$ of Lemma 3.13. The collection of all semilattices of double Ward quasigroups that satisfy 3.131, 3.133 and 3.134 is denoted by $\mathbb{SLDWQ}$.

**Lemma 3.14.** *Suppose that* $(S, *)$ *is a semigroup semilattice* $Y$ *of double Ward quasigroups* $\left(DW_\alpha, */_{DW_\alpha}, e_\alpha\right)$ $(\alpha \in Y)$, *that* $\left(\{e_\alpha : \alpha \in Y\}, *\right) \cong Y$ *and that* $(S, *)$ *satisfies* 3.133 *and* 3.134. *Define* $x_\alpha \bullet y_\beta = \left(e_{\alpha\beta} * x_\alpha\right) * \left(e_{\alpha\beta} * y_\beta\right)$. *Then*

(3.141) $E(S, \bullet) = \left(\{e_\alpha : \alpha \in Y\}, \bullet\right) \cong \left(\{e_\alpha : \alpha \in Y\}, *\right)$ *is a semigroup semilattice,*

(3.142) $(S, \bullet)$ *is a semigroup and a semigroup semilattice* $Y$ *of groups* $\left(DW_\alpha, \bullet/_{DW_\alpha}, e_\alpha\right)$ $(\alpha \in Y)$,

(3.143) *for all* $\alpha, \beta \in Y$ *and all* $x_\alpha \in DW_\alpha$ *and* $y_\beta \in DW_\beta$, $x_\alpha * y_\beta = x_\alpha{}^{-1} \bullet y_\beta{}^{-1}$.

Proof. We have $e_\alpha \bullet e_\beta = \left(e_{\alpha\beta} * e_\alpha\right) * \left(e_{\alpha\beta} * e_\beta\right) = \left[e_{\alpha\beta} * \left(e_\alpha * e_\alpha\right)\right] * \left[e_{\alpha\beta} * \left(e_\beta * e_\beta\right)\right] \overset{3.134}{=} e_\alpha * e_\beta$. By Result 14, $\Psi\left(e_\alpha * e_\beta\right) = \varepsilon = (\Psi\alpha)(\Psi\beta) = \alpha\beta$ and so $e_\alpha * e_\beta = e_{\alpha\beta} = e_\alpha \bullet e_\beta$. This implies that the identity map is an isomorphism from $\left(\{e_\alpha : \alpha \in Y\}, \bullet\right)$ to $\left(\{e_\alpha : \alpha \in Y\}, *\right)$. Furthermore, $e_\alpha \bullet e_\alpha = e_\alpha * e_\alpha = e_\alpha$. Also, $x_\alpha = x_\alpha \bullet x_\alpha = \left(e_\alpha * x_\alpha\right) * \left(e_\alpha * x_\alpha\right)$ implies, by (16), $\left(e_\alpha * x_\alpha\right) * x_\alpha = e_\alpha * x_\alpha$ implies $e_\alpha * x_\alpha = e_\alpha = e_\alpha * e_\alpha$ implies $x_\alpha = e_\alpha$. Hence, $E(S, \bullet) = \{e_\alpha : \alpha \in Y\}$ and so 3.141 is valid.

For each $\left(DW_\alpha, \bullet/_{DW_\alpha}, e_\alpha\right)$ the product is $x_\alpha \bullet y_\alpha = (e_\alpha * x_\alpha) * (e_\alpha * y_\alpha)$ and so, by Result 6, each $\left(DW_\alpha, \bullet/_{DW_\alpha}, e_\alpha\right)$ is a group. By 3.133, $(S,\bullet)$ is a semigroup, so it is a semigroup union of groups whose idempotents commute and, by Lemma 3.1, 3.142 is valid. Finally, by Result 6, $x_\alpha^{-1} = e_\alpha * x_\alpha$ in $(S,\bullet)$. Then, by 3.134, $x_\alpha^{-1} \bullet y_\beta^{-1} = (e_\alpha * x_\alpha) \bullet (e_\beta * y_\beta) =$
$\left[e_{\alpha\beta} * (e_\alpha * x_\alpha)\right] * \left[e_{\alpha\beta} * (e_\beta * y_\beta)\right] \overset{3.134}{=} x_\alpha * y_\beta$. ■

**Definition.** $SLG(S,*)$ denotes $(S,\bullet)$ of Lemma 3.14.

**Theorem 3.15.** *For all* $(S,\bullet) \in \mathbb{SLG}$, $SLG\left(SLDWQ(S,\bullet)\right) = (S,\bullet)$ *and for all* $(S,*) \in \mathbb{SLDWQ}$, $SLDWQ\left(SLG(S,*)\right) = (S,*)$.

Proof. The product in $SLDWQ(S,\bullet)$ is $x * y = x^{-1} \bullet y^{-1}$. The product in $SLG\left(SLDWQ(S,\bullet)\right)$ is
$x_\alpha \otimes y_\beta = \left(e_{\alpha\beta} * x_\alpha\right) * \left(e_{\alpha\beta} * y_\beta\right) = \left(e_{\alpha\beta}^{-1} \bullet x_\alpha^{-1}\right)^{-1} \bullet \left(e_{\alpha\beta}^{-1} \bullet y_\beta^{-1}\right)^{-1} = x_\alpha \bullet \left[e_{\alpha\beta} \bullet \left(y_\beta \bullet e_{\alpha\beta}\right)\right] = \left(x_\alpha \bullet y_\beta\right) \bullet e_{\alpha\beta} = x_\alpha \bullet y_\beta$.
Hence, $SLG\left(SLDWQ(S,\bullet)\right) = (S,\bullet)$.
The product in $SLG(S,*)$ is $x_\alpha \bullet y_\beta = \left(e_{\alpha\beta} * x_\alpha\right) * \left(e_{\alpha\beta} * y_\beta\right)$. The product in $SLDWQ\left(SLG(S,*)\right)$ is
$x \oplus y = x^{-1} \bullet y^{-1} = \left[e_{\alpha\beta} * (e_\alpha * x_\alpha)\right] * \left[e_{\alpha\beta} * \left(e_\beta * y_\beta\right)\right] \overset{8.134}{=} x_\alpha * y_\beta$ and so $SLWD\left(SLG(S,*)\right) = (S,*)$. ■

**Corollary 3.16.** *There is a one-to-one correspondence between semigroup semilattices of groups* ($\mathbb{SLG}$) *and groupoids that are a semigroup semilattice* $Y$ *of double Ward quasigroups and that satisfy* 3.131, 3.133 *and* 3.134 ($\mathbb{SLDWQ}$).

Note that since $\mathbb{SLG}$ is in one-one correspondence with $\mathbb{SLWQ}$, $\mathbb{SLWDQ}$ and $\mathbb{SLDWQ}$, $\mathbb{SLWQ}$ and $\mathbb{SLWDQ}$ are in one-one correspondence with each other, as are $\mathbb{SLDWQ}$ and $\mathbb{SLWQ}$. The next results give the explicit forms of these bijective mappings.

**Theorem 3.17.** $\overleftarrow{\mathbb{SLG}} = \mathbb{SLG}$

Proof. The dual groupoid of a semigroup union of groups with commuting idempotents is a semigroup union of groups with commuting idempotents. As previously noted, the required result then follows from Lemma 3.1. ■

**Theorem 3.18.** $\overleftarrow{\mathbb{SLWQ}} = \mathbb{SLWDQ}$.

Proof. If $(S,*)\in\mathbb{SLWQ}$ then $x*y=x\bullet y^{-1}$ for some $(S,\bullet)\in\mathbb{SLG}$. So, if $\left(T,\overleftarrow{\circ}\right)\in\overleftarrow{\mathbb{SLWQ}}$ then, using 3.17, $x\overleftarrow{\circ}y=y\circ x=y\bullet x^{-1}=x^{-1}\overleftarrow{\bullet}y$ for some $\left(T,\overleftarrow{\bullet}\right)\in\mathbb{SLG}$. As in the proof of Lemma 3.8, $\left(T,\overleftarrow{\circ}\right)$ satisfies 3.81, 3.82, 3.83 and 3.84. Therefore, $\left(T,\overleftarrow{\circ}\right)\in\mathbb{SLWDQ}$. Hence, $\overleftarrow{\mathbb{SLWQ}}\subseteq$ $\mathbb{SLWDQ}$.

If $(S,*)\in\mathbb{SLWDQ}$ then $x*y=x^{-1}\bullet y$ for some $(S,\bullet)\in\mathbb{SLG}$. So, using 3.17, $x\overleftarrow{*}y=y*x=y^{-1}\bullet x=x\overleftarrow{\bullet}y^{-1}$ for some $\left(S,\overleftarrow{\bullet}\right)\in\mathbb{SLG}$. Therefore, as in the proof of Lemma 3.2, $\left(S,\overleftarrow{*}\right)$ satisfies 3.24. Hence, $\left(S,\overleftarrow{*}\right)\in\mathbb{SLWQ}$ and $(S,*)\in\overleftarrow{\mathbb{SLWQ}}$. So, $\mathbb{SLWDQ}\subseteq\overleftarrow{\mathbb{SLWQ}}$. ■

**Corollary 3.19.** $(S,*)\in\mathbb{SLWDQ}$ *if and only if* $(S,*)$ *is a semilattice* $Y$ *of Ward dual quasigroups* $\left(WD_\alpha,*/_{WD_\alpha},e_\alpha\right)(\alpha\in Y)$ *and satisfies the identity* $(x*y)*(z*w)=y*\left[\left(z^{-1}*x^{-1}\right)*w\right]$.

**Theorem 3.20.** $\mathbb{SLDWQ}=\overleftarrow{\mathbb{SLDWQ}}$

Proof. If $(S,*)\in\overleftarrow{\mathbb{SLDWQ}}$ then $x*y=y\overleftarrow{*}x$ for some $\left(S,\overleftarrow{*}\right)\in\mathbb{SLDWQ}$. So, using 3.13 and 3.17, $x*y=y\overleftarrow{*}x=y^{-1}\bullet x^{-1}=x^{-1}\overleftarrow{\bullet}y^{-1}$ for some $\left(S,\overleftarrow{\bullet}\right)\in\mathbb{SLG}$. Therefore, by the proof of Lemma 3.13, $(S,*)\in\mathbb{SLDWQ}$. Hence, $\overleftarrow{\mathbb{SLDWQ}}\subseteq\mathbb{SLDWQ}$ and $\mathbb{SLDWQ}=\overleftarrow{\overleftarrow{\mathbb{SLDWQ}}}\subseteq\overleftarrow{\mathbb{SLDWQ}}$. ■

**Theorem 3.21.** $\mathbb{SLDWQ}$ *and* $\mathbb{SLWQ}$ *are in one-one correspondence.*

Proof. For $(S,*)\in\mathbb{SLDWQ}$ we define $SLWQ(S,*)=(S,\circ)$, where $x_\alpha\circ y_\beta=\left(e_{\alpha\beta}*x_\alpha\right)*y_\beta$. If $(S,\otimes)\in\mathbb{SLWQ}$ we define $SLDWQ(S,\otimes)=(S,\oplus)$, where $x_\alpha\oplus y_\beta=\left(e_{\alpha\beta}\otimes x_\alpha\right)\otimes y_\beta$. Note that, since $(S,*)\in\mathbb{SLDWQ}$, $x*y=x^{-1}\bullet y^{-1}$ for some $(S,\bullet)\in\mathbb{SLG}$. Therefore, $x_\alpha\circ y_\beta=\left(e_{\alpha\beta}*x_\alpha\right)*y_\beta$ $=\left(e_{\alpha\beta}{}^{-1}\bullet x_\alpha{}^{-1}\right)^{-1}\bullet y_\beta{}^{-1}=x_\alpha\bullet e_{\alpha\beta}\bullet y_\beta{}^{-1}=x_\alpha\bullet e_\alpha\bullet e_\beta\bullet y_\beta{}^{-1}=x_\alpha\bullet y_\beta{}^{-1}$. By Lemma 3.2, $SLWQ(S,*)=(S,\circ)\in\mathbb{SLWQ}$. Therefore, $SLDWQ(S,\circ)=(S,\oplus)$, where $x_\alpha\oplus y_\beta=\left(e_{\alpha\beta}\circ x_\alpha\right)\circ y_\beta=\left(e_{\alpha\beta}\bullet x_\alpha{}^{-1}\right)\bullet y_\beta{}^{-1}=x_\alpha{}^{-1}\bullet y_\beta{}^{-1}$ and so $(S,\oplus)\in\mathbb{SLDWQ}$. So, $SLDWQ:\mathbb{SLWQ}\to\mathbb{SLDWQ}$ and $SLWQ:\mathbb{SLDWQ}\to\mathbb{SLWQ}$.

Then, the product in $SLDWQ(SLWQ(S,*))$ is $x_\alpha \oplus y_\beta = (e_{\alpha\beta} \circ x_\alpha) \circ y_\beta = [e_{\alpha\beta} * (e_{\alpha\beta} \circ x_\alpha)] * y_\beta =$
$= [e_{\alpha\beta} * (e_{\alpha\beta} * x_\alpha)] * y_\beta = \left[e_{\alpha\beta}{}^{-1} \bullet (e_{\alpha\beta}{}^{-1} \bullet x_\alpha{}^{-1})^{-1}\right]^{-1} \bullet y_\beta{}^{-1} = \left[(e_{\alpha\beta}{}^{-1} \bullet x_\alpha{}^{-1})^{-1}\right]^{-1} \bullet y_\beta{}^{-1} = (e_{\alpha\beta}{}^{-1} \bullet x_\alpha{}^{-1}) \bullet y_\beta{}^{-1} =$
$= x_\alpha{}^{-1} \bullet y_\beta{}^{-1} = x_\alpha * y_\beta$. Therefore, $SLDWQ(SLWQ(S,*)) = (S,*)$. Similarly, $SLWQ(SLDWQ(S,\otimes)) = (S,\otimes)$. ■

**Questions.** Suppose that $(S,*)$ is a semigroup semilattice $Y$ of double Ward quasigroups $(DW_\alpha, */_{DW_\alpha}, e_\alpha)$ $(\alpha \in Y)$ and that $(S,*)$ satisfies 3.133 and 3.134. Then, is $(\{e_\alpha : \alpha \in Y\}, *) \cong Y$ .

Can groupoids in $\mathbb{SLDWQ}$ be described by a single identity, in place of 3.133 and 3.134?

Is there a structure theorem for groupoids in $\mathbb{SLWQ}$, $\mathbb{SLWDQ}$ and $\mathbb{SLDWQ}$ analogous to the structure theorem for semigroups that are semigroup semilattices of groups [3, Theorem 4.11]?

A remaining area for investigation is right, left and double division on rectangular groups, completely simple semigroups and semigroup unions of groups, where $x^{-1}$ is the inverse of $x$ in the group to which it belongs.

## References

[1] Cardoso, J. M. and da Silva, C. P.: *On Ward quasigroups*, An. Stiint. Univ. “Al. I. Cuza” Iasi Sect. I a Mat. (N.S.) **24** (1978), no.2, 231-233.

[2] Chatterjea, S. K.: *On Ward quasigroups*, Pure Math. Manuscript **6** (1987), 31-34.

[3] Clifford, A. H. and Preston, G. B.: *The algebraic theory of semigroups*, vol.1. In: Mathematical Survey, vol.7. American Mathematical Society, Providence (1961).

[4] Da Silva, C. P.: *On a theorem of Lagrange for Ward quasigroups* (Portugese), Rev. Columbiana Ma. **12** (1978), no.3-4, 91-96.

[5] Fiala, N. C.: *Double Ward quasigroups*, Quasigroups and Related Systems **15** (2007), 261-262.

[6] Furstenberg, H.: *The inverse operation in groups*, Proceedings of the American Mathematical Society **6**, no.6 (Dec. 1955), 991-997.

[7] Higman, G. and Neumann, B.H.: *Groups as groupoids with one law*, Publ. Math. Debrecen **2** (1952), 215-221.

[8] Johnson, K. W. and Vojtechovsky, P.: *Right division in groups, Dedikind-Frobenius group matrices and Ward quasigroups*, Abh. Math. Sem. Univ. Hamburg **75** (2006), 121-36

[9] Monzo, R. A. R.: *The ternary operations of groupoids*, arXiv:1510.07955 [math. RA].

[10] Morgado, J.: *Definicao,de quasigrupo subtractive por um unico axioma*, Gaz. Mat. (Lisboa) **92-93** (1963), 17-18.

[11] Polonijo, M.: *A note on Ward quasigroups*, An. Stiint. Univ. "Al. I. Cuza" Iasi Sect. I a Mat. (N.S.) 32 (1986), no.2, 5-10.

[12] Polonijo, M.: *Transitive Groupoids*, Portugaliae Mathematica, **50** (1993), no. 1, 63-74.

[13] Rabinow, D. G.: *Independent sets of postulates for abelian groups and fields in terms of the inverse operations*, American Journal of Mathematics **59**, no.1 (Jan. 1937), 211-224.

[14] Ward, M.: *Postulates for the inverse operations in a group*, Transactions of the American Mathematical Society, **32, no.3,** (July 1930), 520-526.

[15] Whittaker, J. V.: *On the postulates defining a group*, Amer. Math. Monthly **62** (1955), 636-640.

R. A. R. Monzo
Flat 10 Albert Mansions
Crouch Hill, London N89RE
United Kingdom
bobmonzo@talktalk.net